\documentclass[10pt]{amsart}
\usepackage{amssymb}

\newcommand{\expan}{{\mathrm{ex}}}

\newcommand{\csmf}{compact smooth manifold}

\theoremstyle{definition}
\newtheorem{thm}{Theorem}[section]
\newtheorem{lem}[thm]{Lemma}
\newtheorem{prp}[thm]{Proposition}
\newtheorem{dfn}[thm]{Definition}
\newtheorem{cor}[thm]{Corollary}
\newtheorem{cnj}[thm]{Conjecture}

\newcommand{\beq}{\begin{equation}}
\newcommand{\eeq}{\end{equation}}
\newcommand{\beqr}{\begin{eqnarray*}}
\newcommand{\eeqr}{\end{eqnarray*}}
\newcommand{\bal}{\begin{align*}}
\newcommand{\eal}{\end{align*}}
\newcommand{\bei}{\begin{itemize}}
\newcommand{\eei}{\end{itemize}}

\newcommand{\gm}{\gamma}
\newcommand{\dt}{\delta}
\newcommand{\ep}{\varepsilon}
\newcommand{\zt}{\zeta}
\newcommand{\et}{\eta}

\newcommand{\ld}{\lambda}

\newcommand{\ph}{\varphi}
\newcommand{\ps}{\psi}

\newcommand{\ta}{\tau}

\newcommand{\Gm}{\Gamma}

\newcommand{\Z}{{\mathbf{Z}}}
\newcommand{\R}{{\mathbf{R}}}
\newcommand{\C}{{\mathbf{C}}}

\newcommand{\ev}{{\mathrm{ev}}}
\newcommand{\sint}{{\mathrm{int}}}

\newcommand{\spec}{{\mathrm{sp}}}

\newcommand{\diag}{{\mathrm{diag}}}

\newcommand{\rank}{{\mathrm{rank}}}

\newcommand{\dirlim}{\displaystyle \lim_{\longrightarrow}}

\newcommand{\andeqn}{\,\,\,\,\,\, {\mbox{and}} \,\,\,\,\,\,}
\newcommand{\QED}{\rule{0.4em}{2ex}}

\newcommand{\ca}{C*-algebra}
\newcommand{\ct}{continuous}
\newcommand{\pj}{projection}
\newcommand{\mf}{manifold}

\newcommand{\nbhd}{neighborhood}

\newcommand{\hm}{homomorphism}

\newcommand{\ifo}{if and only if}
\newcommand{\rsha}{recursive subhomogeneous algebra}

\newcommand{\rshd}{recursive subhomogeneous decomposition}

\newcommand{\tdim}{topological dimension}

\newcommand{\mh}{minimal homeomorphism}

\newcommand{\cms}{compact metric space}

\newcommand{\hsa}{hereditary subalgebra}

\newcommand{\scn}{strong covering number}

\newcommand{\rsz}[1]{\raisebox{0ex}[0.8ex][0.8ex]{$#1$}}
\newcommand{\subrsz}[1]{\raisebox{0ex}[0.8ex][0.8ex]{${\scriptstyle{#1}}$}}

\title[C*-algebras of minimal
diffeomorphisms]{Direct limit decomposition for C*-algebras
of minimal diffeomorphisms}

\author{Qing  Lin}

\address{Ericsson Canada Inc.,
   18th Floor, 1140 West Pender St.,
   Vancouver BC V6E 4G1, Canada.}

\email[]{qing.lin@ericsson.ca}

\author{N.\  Christopher Phillips}

\address{Department of Mathematics, University  of Oregon,
       Eugene OR 97403-1222, USA.}

\email[]{ncp@darkwing.uoregon.edu}

\thanks{Research of the second author
      partially supported by NSF grants DMS 9400904 and DMS 9706850.}

\begin{document}

\maketitle

This article outlines the proof that the crossed product
$C^* (\Z, M, h)$ of a \csmf\  $M$ by a minimal diffeomorphism
$h : M \to M$ is isomorphic to a direct limit of
subhomogeneous \ca s belonging to a tractable class.
This result is motivated by the Elliott classification program
for simple nuclear \ca s \cite{El2}, and the observation that
the known classification theorems in the stably finite case
mostly apply to certain kinds of direct limits of subhomogeneous \ca s,
or at least to \ca s with related structural conditions.
(See Section~1.)
This theorem is a generalization, in a sense, of direct limit
decompositions for
crossed products by \mh s of the Cantor set (Section~2
of~\cite{Pt2}), for the irrational rotation algebras (\cite{EE}),
and for some higher dimensional noncommutative toruses
(\cite{EL1}, \cite{EL2}, \cite{LQ1}, and~\cite{Bc}).
(In \cite{Pt2}, only a local approximation result is stated,
but the \ca s involved are semiprojective.)
Our theorem is not a generalization in the strict sense for several
reasons; see the discussion in Section~1.

There are four sections.
In the first, we state the theorem and discuss some consequences
and expected consequences.
In the second section, we describe the basic construction in our proof,
a modified Rokhlin tower, and show how \rsha s appear naturally
in our context.
The third section describes how to prove local approximation by
\rsha s, a weak form of the main theorem.
In Section~4, we give an outline of how to use the methods of
Section~3 to obtain the direct limit decomposition.

This paper is based on a talk given by the second author at the
US--Japan Seminar on Operator Algebras and Applications
(Fukuoka, June 1999), which roughly covered Sections~2 and~3,
and on a talk given by the second author at the 28th Canadian
Annual Symposium on Operator Algebras (Toronto, June 2000),
which roughly covered Sections~1 and~2.
At the time of the first talk, only the local approximation result
described in Section~3 had been proved.
We refer to the earlier survey paper \cite{LP1} for earlier
parts of the story; this paper reports the success of the project
described in Section~6 there.

The first author would like to thank George Elliott, John Phillips, and 
Ian Putnam for funding him at the University of Victoria where
some of this work was carried out.
He would particularly like to acknowledge his great 
gratitude to Ian Putnam for many interesting discussions.
The second author would like to thank Larry Brown,
Marius D\v{a}d\v{a}rlat, George Elliott,
and Ian Putnam for useful
discussions and email correspondence.
Some of the work reported here was
carried out during a sabbatical year at Purdue University,
and he would like to thank that institution for its hospitality.

\section{The main theorem, consequences, and conjectured consequences}

The main theorem is as follows.
Undefined terminology is discussed after the statement.

\begin{thm}\label{main}
Let $M$ be a connected
\csmf\  with $\dim (M) = d > 0$, and let $h \colon M \to M$
be a minimal diffeomorphism.
Then there exists an increasing sequence
\[
A_0 \subset A_1 \subset A_2 \subset \cdots \subset C^* (\Z, M, h)
\]
of C*-subalgebras of $C^* (\Z, M, h)$ such that
\[
\overline{\textstyle{\bigcup_{n = 0}^{\infty} A_n}} = C^* (\Z, M, h)
\]
and such that each $A_n$ has a separable \rshd\  with \tdim\  at most
$d$ and \scn\  at most $d (d + 2)$.
\end{thm}

A \rsha\  (a \ca\  with a \rshd)
is a particularly tractable kind of subhomogeneous \ca.
See \cite{Ph6}, \cite{Ph7}, and \cite{LP1}, and also see the consequences
below.
We will explain in Section~2 how \rsha s arise, and we will recall
(informally) the definition there (after Theorem~\ref{AYS-rshd}).
A finite direct sum
\[
\bigoplus_{k = 0}^l C {\textstyle{ \left( X_k, M_{n (k)} \right) }}
\]
of (trivial) homogeneous \ca s is a special case of a \rsha,
and the \tdim\  is simply $\max_{0 \leq k \leq l} \dim (X_k)$.
(Dimension is taken to be covering dimension;
see Definition 1.6.7 of \cite{En}.)
The condition in the theorem that $A_n$ have \tdim\  at most $d$
for all $n$ thus ensures that the resulting direct limit decomposition
$C^* (\Z, M, h) \cong \dirlim A_n$ has no dimension growth.

In general, it is not possible to find
a representation as a direct limit
(with no dimension growth) of direct sums of
corners of trivial homogeneous \ca s.
A simple direct limit of this sort must even be approximately
divisible in the sense of \cite{BKR}, by Theorem~2.1 of
\cite{EGL1}.
However, a crossed product by a minimal diffeomorphism may have
no nontrivial \pj s
(Corollary~3 and Example~4 of Section~5 of \cite{Cn}).

We will not define the  \scn\  here, although some discussion will
be given after Theorem~\ref{EL-main}.
We have included it in the conclusion because
the proof of Theorem~\ref{EL-main} suggests that a bound on
the \scn\  might be necessary for some classification results.

The requirement that we have a diffeomorphism of a manifold is
connected with the appearance of a condition on the \scn\  in the
hypotheses of Theorem~\ref{EL-main}.
This also will be discussed after that theorem.
We certainly expect that the theorem will be true for minimal
homeomorphisms of finite dimensional compact metric spaces
(even, presumably,
compact metric spaces with infinite covering dimension).

We point out here that our theorem does not directly imply the
Elliott-Evans direct limit representation for the
irrational rotation algebras \cite{EE}.
Our theorem gives a representation of an irrational rotation algebra
as a direct limit of \rsha s with \tdim\  at most $1$,
while the Elliott-Evans theorem gives a representation
as a direct limit of direct sums of homogeneous \ca s with
\tdim\  at most $1$ (in fact, circle algebras).
We do not recover the results of~\cite{EL1}, \cite{EL2}, and~\cite{LQ1}
(for certain higher dimensional noncommutative toruses),
not only because the algebras in our direct system are more
complicated but also because not all the algebras considered there
are even crossed products by diffeomorphisms.
We also do not recover the direct limit decomposition
for crossed products by minimal homeomorphisms of the Cantor set
(see Section~2 of~\cite{Pt2} for the local approximation result),
because the Cantor set is not a \mf.
(Our methods do specialize to this case, but that would be silly, since
our argument is much more complicated.)

Theorem~\ref{main} has the following consequences for crossed products
by minimal diffeomorphisms.
These consequences all hold for an arbitrary simple unital direct limit
of \rsha s, assuming no dimension growth and that the maps of
the system are unital and injective.
(Most don't require the full strength of these hypotheses,
but all require some restriction on dimension growth.
None require any hypotheses on the \scn.)
The proofs are in \cite{Ph7},
and the statements can be found in Section~4 of \cite{LP1}
(except for the last one, which is actually a consequence of stable
rank one).
In all of these, $M$ is a connected
\csmf\  with $\dim (M) > 0$, and $h \colon M \to M$
is a minimal diffeomorphism.

\begin{cor}\label{tsr}
(Theorem~3.6 of \cite{Ph7}.)
The algebra $C^* (\Z, M, h)$ has stable rank one in the sense of
\cite{Rf}.
That is, the invertible group ${\rm{inv}} ( C^* (\Z, M, h) )$ is dense
in $C^* (\Z, M, h)$.
\end{cor}

\begin{cor}\label{canc}
(Theorem~2.2 of \cite{Ph7}.)
The \pj s in
\[
M_{\infty} (C^* (\Z, M, h))
            = \bigcup_{n = 1}^{\infty} M_n (C^* (\Z, M, h))
\]
satisfy cancellation.
That is, if $e, \, p, \, q \in M_{\infty} (C^* (\Z, M, h))$ are \pj s,
and if $p \oplus e \sim q \oplus e$, then $p \sim q$.
\end{cor}

\begin{cor}\label{FCQ}
(Theorem~2.3 of \cite{Ph7}.)
The algebra $C^* (\Z, M, h)$ satisfies
Blackadar's Second Fundamental Comparability Question
(\cite{Bl3}, 1.3.1).
That is, if $p, \, q \in M_{\infty} (C^* (\Z, M, h))$ are \pj s,
and if $\ta (p) < \ta (q)$ for every normalized trace $\ta$ on
$C^* (\Z, M, h)$, then $p \precsim q$.
\end{cor}

\begin{cor}\label{unperf}
(Theorem~2.4 of \cite{Ph7}.)
The group $K_0 (C^* (\Z, M, h))$ is unperforated for the strict order.
That is, if $\et \in K_0 (C^* (\Z, M, h))$ and if there is $n > 0$
such that $n \et > 0$, then $\et > 0$.
\end{cor}

(In the simple case, this is the same as saying that
$K_0 ( C^* (\Z, M, h) )$
is weakly unperforated in the sense of 2.1 of \cite{El1}.)

\begin{cor}\label{UmodU0}
(Theorem~2.1 of \cite{Ph7}.)
The canonical map
\[
U (C^* (\Z, M, h)) / U_0 (C^* (\Z, M, h)) \to K_1 (C^* (\Z, M, h))
\]
is an isomorphism.
\end{cor}

A small part of these results could already be obtained using the
weaker (and much simpler) methods described in Sections~1 and~5
of \cite{LP1}.
For example, it had already been shown
that the order on $K_0 (C^* (\Z, M, h))$ is
determined by traces (a weak form of Corollary~\ref{FCQ}),
and hence that
$K_0 (C^* (\Z, M, h))$ is unperforated for the strict order
(Corollary~\ref{unperf}).
Also, surjectivity in Corollary~\ref{UmodU0} (but not
injectivity) was known.

The criterion in \cite{BDR}, for when a simple direct limit of
direct sums of trivial homogeneous \ca s with slow
dimension growth
has real rank zero, is known to fail for simple direct limits
of \rsha s with no dimension growth.
(Indeed, it even fails for crossed products by minimal diffeomorphisms;
see Example~5.7 of \cite{LP1}.)
Nevertheless, it appears likely that a suitable strengthening of the
condition will be equivalent to real rank zero for such direct limits,
and that the proof will not be difficult.
Specializing (for simplicity) to the case of a unique trace,
we obtain the following, which we state as a conjecture.

\begin{cnj}\label{RR0}
Let $M$ be a connected
\csmf\  with $\dim (M) > 0$, and let $h \colon M \to M$
be a uniquely ergodic minimal diffeomorphism.
Let $\ta \colon C^* (\Z, M, h) \to \C$ be the trace induced by
the unique invariant probability measure.
Then $C^* (\Z, M, h)$ has real rank zero (\cite{BP})
\ifo\  $\ta_* ( K_0 ( C^* (\Z, M, h)) )$ is dense in $\R$.
\end{cnj}

For methods for computing the ranges of traces on the K-theory of
crossed products by $\Z$, we refer to \cite{Ex2}.

It might not be terribly difficult to prove that if a simple \ca\  $A$
is a direct limit of a system of \rsha s with no dimension growth,
and possibly also assuming that the maps of the system are unital and
injective, then real rank zero implies tracial rank zero in the sense
of H.\  Lin \cite{Ln14}.
If so, then the following result of H.\  Lin (Theorem~3.9 of~\cite{Ln15})
implies classifiability:

\begin{thm}\label{class}
Let $A$ and $B$ be separable simple unital \ca s
with tracial rank zero in the sense of \cite{Ln14}.
Suppose that each has local approximation by subalgebras with
bounded dimensions of irreducible representations.
That is, for every finite subset $F \subset A$ and every $\ep > 0$,
there is a C*-subalgebra $D \subset A$ and an integer $N$
such that every
element of $F$ is within $\ep$ of an element of $D$ and every irreducible
representation of $D$ has dimension at most $N$; and similarly for $B$.
Then
\[
{\textstyle{
\left( K_0 (A), \, K_0 (A)_+, \, [1_A], \, K_1 (A) \right)
\cong \left( K_0 (B), \, K_0 (B)_+, \, [1_B], \, K_1 (B) \right)    }}
\]
implies $A \cong B$.
\end{thm}

In particular, one would have a proof of the following conjecture:

\begin{cnj}\label{cpcl}
Let $M$ be a connected
\csmf\  with $\dim (M) > 0$, and let $h \colon M \to M$
be a uniquely ergodic minimal diffeomorphism.
Let $\ta \colon C^* (\Z, M, h) \to \C$ be the trace induced by
the unique invariant probability measure, and assume that
$\ta_* ( K_0 ( C^* (\Z, M, h)) )$ is dense in $\R$.
Then the crossed product \ca\  $C^* (\Z, M, h)$
is classifiable.
\end{cnj}

We will not give a precise definition of ``classifiable'' here.

We note that H.\  Lin's classification theorem has no hypotheses
involving slow dimension growth, and does not even require
a direct limit representation; only local approximation is needed,
and the condition on the approximating algebras is weak.
(Indeed, H.\  Lin has other classification theorems which don't
even require local approximation, but do require further restrictions
on the K-theory.)
However, at least with our current state of knowledge, the
direct limit representation in Theorem~\ref{main}, including
the no dimension growth condition, seems to be needed to verify
the other hypotheses of Theorem~\ref{class}.
For example, simple direct limits that don't have slow dimension growth
need not even have stable rank one \cite{Vl}.

Since $C^* (\Z, M, h)$ always has stable rank one,
if it doesn't have real rank zero then it has real rank one.
However, most of the currently known general classification theorems
apply only to algebras with many projections,
and those that don't are much too restrictive in other ways
(such as assuming trivial K-theory).
In particular, the \ca s covered by \cite{Gn} and \cite{EGL2}
are approximately divisible
(as discussed above), and the theorems of H.\  Lin
(see \cite{Ln12} and \cite{Ln13})
require a finite value of the tracial rank,
the definition of which again requires
the existence of many nontrivial \pj s.
However, as mentioned above, the example of Connes shows that
$C^* (\Z, M, h)$ may have no nontrivial \pj s.
There is a classification theorem \cite{JS} for a special class of
direct limits which includes simple \ca s with no nontrivial \pj s,
but the building blocks there are much more special than those
appearing in our theorem.

We are hopeful that the approach of \cite{Gn} and \cite{EGL2},
which now covers
simple direct limits, with no dimension growth,
of direct sums of homogeneous \ca s (actually, a slightly larger class),
can be generalized to cover simple direct limits,
with no dimension growth, of \rsha s, possibly with the added
restriction of no growth of the \scn.
One reason for optimism (as well as for the belief that conditions
on the \scn\  might be necessary) is the successful generalization of
exponential length results from the case of trivial homogeneous \ca s
to \rsha s; see Theorem~\ref{EL-main} below.
The related results for the trivial homogeneous case (see 
Theorems~3.3 and~4.5 of \cite{Ph1.5}) depended heavily on the existence
of many \pj s, but in the proof of Theorem~\ref{EL-main} we had
to learn to handle situations with no nontrivial \pj s at all.
However, we do not know whether Theorem~\ref{EL-main} is
even true without the condition on the \scn.
(See the discussion after the statement of that theorem.)
We included the bound on the \scn\  in Theorem~\ref{main}
because of the possibility that it might be necessary for our
suggested approach to proving a classification result, or
perhaps even for a classification result to hold.

In any case,
a generalization of the methods of \cite{Gn} and \cite{EGL2}
is likely to be very difficult.
Possibly the situation will be improved by a generalization of
H.\  Lin's methods that is strong enough to apply to simple \ca s
which contain no nontrivial \pj s.

\section{Modified Rokhlin towers}

Throughout this section, $M$ is a \cms\  and $h \colon M \to M$ is a \mh.
(The requirement that $M$ be a \mf\  will not be needed until
the next section.)
We let $u$ denote the implementing unitary in $C^* (\Z, M, h)$,
so that $u f u^* = f \circ h^{-1}$ for $f \in C (M)$.

We start with a definition.

\begin{dfn}\label{FirstRet}
Let $Y \subset M$, and let $x \in Y$.
The {\emph{first return time}} $\ld_Y (x)$
(or $\ld (x)$ if $Y$ is understood)
of $x$ to $Y$
is the smallest integer $n \geq 1$ such that $h^n (x) \in Y$.
We set $\ld (x) = \infty$ if no such $n$ exists.
\end{dfn}

The following result is well known in the area, and is easily proved:

\begin{lem}\label{RokhFin}
If $\sint (Y) \neq \varnothing$, then $\sup_{x \in Y} \ld (x) < \infty$.
\end{lem}

Let $Y \subset M$ with $\sint (Y) \neq \varnothing$.
Let $n (0) < n (1) < \cdots < n (l)$ (or, if the dependence on $Y$ must be
made explicit, $n_Y (0) < n_Y (1) < \cdots < n_Y (l_Y)$)
be the distinct values of $\ld (x)$ for $x \in Y$.
The Rokhlin tower based on a subset
$Y \subset M$ with $\sint (Y) \neq \varnothing$ consists of the partition
\[
Y = \coprod_{k = 0}^{l} \{ x \in Y \colon \ld (x) = n (k) \}
\]
of $Y$ (the sets here are the base sets), and the corresponding
partition
\[
M = \coprod_{k = 0}^{l} \coprod_{j = 0}^{n (k) - 1}
    h^j {\textstyle{ \left( {\rsz{
             \{ x \in Y \colon \ld (x) = n (k) \} }} \right) }}
\]
of $M$.
Note that $h$ acts like a cyclic shift except on the top space
\[
h^{n (k) - 1} {\textstyle{ \left( {\rsz{
             \{ x \in Y \colon \ld (x) = n (k) \} }} \right) }}
\]
of each ``tower''
\[
\coprod_{j = 0}^{n (k) - 1}
    h^j {\textstyle{ \left( {\rsz{
             \{ x \in Y \colon \ld (x) = n (k) \} }} \right) }}.
\]

Actually, for our purposes it is more convenient to use the partition
\[
M = \coprod_{k = 0}^{l} \coprod_{j = 1}^{n (k)}
    h^j {\textstyle{ \left( {\rsz{
             \{ x \in Y \colon \ld (x) = n (k) \} }} \right) }}.
\]
Note that
\[
Y = \coprod_{k = 0}^{l} h^{n (k)}
         {\textstyle{ \left( {\rsz{
             \{ x \in Y \colon \ld (x) = n (k) \} }} \right) }},
\]
so that $h$ now acts like a cyclic shift on the towers,
except on $Y$ itself.

We will be interested in arbitrarily small choices for $Y$,
in particular with arbitrarily small diameter and for which the
smallest first return time $n_Y (0)$ is arbitrarily large.
If $M$ is totally disconnected, then we may choose $Y$ to be
both closed and open.
In this case, the sets
\[
Y_k = \{ x \in Y \colon \ld (x) = n (k) \}
\]
are all closed, and
there is a composite \hm\  $\gm_0$ given by
\[
C (M) \longrightarrow \bigoplus_{k = 0}^{l} \bigoplus_{j = 1}^{n (k)}
     C (h^j ( Y_k ) )
       \stackrel{\cong}{\longrightarrow}
   \bigoplus_{k = 0}^{l}  C ( Y_k )^{n (k)},
\]
which is in fact an isomorphism.
The formula is
\[
\gm_0 (f) = \left( {\rsz{ {\textstyle{
   \left( f \circ h |_{Y_0}, \, \dots, \, f \circ h^{n (0)} |_{Y_0} \right),
   \, \dots, \, 
   \left( f \circ h |_{Y_l}, \, \dots, \, f \circ h^{n (l)} |_{Y_l} \right)
    }} }} \right).
\]
See \cite{Pt1} for the exploitation of this idea.

In order to have a C*-algebraically sensible codomain for $\gm_0$,
we must insist that the sets $Y_k$ be closed.
However, the spaces $M$ we are interested in are connected, so
we are forced to choose
\[
Y_k = \overline{ \{ x \in Y \colon \ld (x) = n (k) \} }
\]
instead.
The sets $h^j ( Y_k )$ are no longer disjoint (although they certainly
cover $M$), so our map
\[
\gm_0 \colon C (M) \to \bigoplus_{k = 0}^{l}  C ( Y_k )^{n (k)},
\]
while still injective, is no longer an isomorphism.

Next, define
\[
s_k \in M_{n (k)} \subset C \left( Y_k, M_{n (k)} \right)
\]
by
\[
s_k = \left( \begin{array}{ccccccc}
         0 &      0 & \cdots & \cdots &      0 &      0 &      1 \\
         1 &      0 & \cdots & \cdots &      0 &      0 &      0 \\
         0 &      1 & \cdots & \cdots &      0 &      0 &      0 \\
    \vdots & \vdots & \ddots &        & \vdots & \vdots & \vdots \\
    \vdots & \vdots &        & \ddots & \vdots & \vdots & \vdots \\
         0 &      0 & \cdots & \cdots &      1 &      0 &      0 \\
         0 &      0 & \cdots & \cdots &      0 &      1 &      0
 \end{array} \right),
\]
and define
\[
s = (s_0, s_1, \dots, s_l )
   \in \bigoplus_{k = 0}^{l} C \left( Y_k, M_{n (k)} \right).
\]
Then $s$ is unitary.
Identifying $C ( Y_k )^{n (k)}$ with the diagonal matrices in
$C \left( Y_k, M_{n (k)} \right)$ in the obvious way, one can
check that if $f \in C (M)$ vanishes on $Y$, then
\[
\gm_0 (u f u^*) = \gm_0 (f \circ h^{-1}) = s \gm_0 (f) s^*.
\]
The calculation uses the fact that
\[
Y = \bigcup_{k = 0}^{l} h^{n (k)} ( Y_k),
\]
and in fact our choice to start our towers at $h (Y_k)$ rather than
at $Y_k$ was made to have this formula work correctly when
$f$ vanishes on $Y$ (rather than when $f$ vanishes on $h^{-1} (Y)$).

This relation allows us to extend $\gm_0$ to the following subalgebra
of $C^* (\Z, M, h)$:

\begin{dfn}\label{AY-def}
For any closed subset $Y \subset M$, we define
\[
A (Y)
  = C^* {\textstyle{ \left( {\rsz{ C (M), \, u C_0 (M \setminus Y) }}
                                                 \right)  }}
         \subset C^* (\Z, M, h),
\]
the C*-subalgebra of $C^* (\Z, M, h)$ generated by $C (M)$ and
$u C_0 (M \setminus Y)$.
Here, we identify $C_0 (M \setminus Y)$ in the obvious way
with the subalgebra of $C (M)$ consisting of those functions vanishing
on $Y$.
We use the analogous convention throughout the paper.
\end{dfn}

\begin{prp}\label{AY-hom}
Let $M$ be a \cms, and let $h \colon M \to M$ be a \mh.
Let $Y \subset M$ be closed with $\sint (Y) \neq \varnothing$.
Then there exists a unique \hm\  %
\[
\gm_Y \colon A (Y) \to
     \bigoplus_{k = 0}^l C \left( Y_k, M_{n (k)} \right)
\]
such that if $f \in C (M)$, then
\[
{\textstyle{
\gm_Y ( f )_k = \diag \left( f \circ h |_{Y_k},
    \, f \circ h^2 |_{Y_k}, \, \dots,
    \, f \circ h^{n (k)} |_{Y_k} \right)   }}
\]
and if $f \in C_0 (M \setminus Y)$, then
\[
(\gm_Y (u f))_k = s_k \gm_Y (f)_k.
\]
Moreover, $\gm_Y$ is unital and injective.
\end{prp}

We now introduce a slight twist on these ideas.

\begin{dfn}\label{AYS-def}
Let $Y \subset M$ be closed with $\sint (Y) \neq \varnothing$.
Let $S \subset \sint (Y_0)$ be closed.
Define
\[
\expan (S) =
 {\textstyle{ \left\{ h (S), h^2 (S), \dots, h^{n (0)} (S) \right\} }},
\]
which is a collection of disjoint closed subsets of $M$.
Define $C (M)_{\expan (S)}$ to be the set of all $f \in C (M)$
such that $f$ is constant on $T$ for every $T \in \expan (S)$.
(The constant value is allowed to depend on $T$.)
Define $A (Y, S)$ to be the C*-subalgebra of
$C^* (\Z, M, h)$ given by
\[
{\textstyle{
A (Y, S) = C^* \left( C (M)_{ \expan (S) }, \,
 u \left[ \rsz{ C_0 (M \setminus Y) \cap C (M)_{ \expan (S) } }
                             \right] \right)
 \subset A (Y).   }}
\]
\end{dfn}

As we will see below, the point of this definition is that
(when $\sint (S) \neq \varnothing$) we
can construct useful unitaries in $C^* (\Z, M, h)$ which
commute with $A (Y, S)$.
(See Step~9 in the proof outline in Section~3.)

It is not obvious what the image
\[
\gm_Y (A (Y, S))
   \subset \bigoplus_{k = 0}^l C \left( Y_k, M_{n (k)} \right)
\]
looks like,
and working with it directly threatens to be very complicated.
Fortunately, the essential properties can be abstracted in
a tractable way; the result is what we call a \rsha.
(The definition of a \rsha\  was
in fact invented for exactly this purpose.)
First, we recall the notion of a pullback.

\begin{dfn}\label{A0}
Let $A$ and $B$ be \ca s, and let a third \ca\  $C$ and
\hm s $\ph \colon A \to C$ and $\ps \colon B \to C$ be given.
The {\emph{pullback}} (also called fibered product or restricted
direct sum) is
\[
A \oplus_C B =  A \oplus_{C, \ph, \ps} B =
\{ (a, b) \in A \oplus B \colon \ph (a) = \ps (b) \}.
\]
If the maps $\ph$ and $\ps$ are understood,
we will write  $A \oplus_C B$.
\end{dfn}

\begin{thm}\label{AYS-rshd}
Let $M$ be a \cms, and let $h \colon M \to M$ be a \mh.
Let $Y \subset M$ be closed with $\sint (Y) \neq \varnothing$.
Let $S \subset \sint (Y_0)$ be closed.
Then there exist closed subsets
\[
Y_k^{(0)} \subset \partial Y_k \subset Y_k
\]
for $1 \leq k \leq l$,
and \hm s $\ph_k$ and $\ps_k$ (with $\ps_k$ being just the
restriction map) such that the image $\gm_Y (A (Y, S))$ is
equal to the subalgebra
\begin{align*}
\left[ \cdots \rule{0em}{3ex} \left[ \rule{0em}{2.9ex} \left[
  {\textstyle{  C \left( Y_0, M_{n (0)} \right)_S }}
              \right.  \right.  \right. \!\!\!  & \left. \left. \left.
      \oplus_{ C \left( \subrsz{ Y_1^{(0)}, M_{n (1)} }
                                  \right), \ph_1, \ps_1 }
        {\textstyle{ C \left( Y_1, M_{n (1)} \right) }} \right]
             \rule{0em}{2.9ex}  \right. \rule{0em}{3ex} \right.  \\
    & \left. \left.
    \oplus_{ C \left( \subrsz{ Y_2^{(0)}, M_{n (2)} }
                                 \right), \ph_2, \ps_2 }
               {\textstyle{ C \left( Y_2, M_{n (2)} \right) }}
      \rule{0em}{2.9ex} \right] \cdots  \rule{0em}{3ex} \right] \\
    &  \left. \left. \oplus_{ C \left( \subrsz{ Y_l^{(0)}, M_{n (l)} }
                                 \right), \ph_l, \ps_l }
                {\textstyle{ C \left( Y_l, M_{n (l)} \right) }}
                        \right.  \right.
\end{align*}
of $\bigoplus_{k = 0}^l C \left( Y_k, M_{n (k)} \right)$.
Here, by analogy with Definition~\ref{AYS-def}, we set
\[
C \left(Y_0, M_{n (0) } \right)_S
  = \left\{ f \in C \left( Y_0, M_{n (0) } \right) \colon
     {\mbox{$f$ is constant on $S$}} \right\}.
\]
\end{thm}

A \ca\  given as an iterated pullback as in the conclusion of this
theorem, in which the algebras have the form
$C \left( X_k, M_{n (k} \right)$,
the maps $\ph_k$ are unital,
and the maps $\ps_k$ are unital and surjective,
is called a {\emph{\rsha}}.
We refer to Section~2 of \cite{LP1} for a more careful definition,
for some useful associated terminology, and examples;
to Section~3 of \cite{LP1} for a discussion of the proof of
Theorem~\ref{AYS-rshd} (in the case $S = \varnothing$);
and to Section~4 of \cite{LP1} for a discussion of why the concept
of a \rshd\  is useful and what can be done with it.
We recall here that the \tdim\  is the largest dimension $\dim (X_k)$.
Unfortunately, it depends on the particular decomposition; see
Example~2.9 of \cite{LP1}.
We will always have a decomposition in mind, usually coming from
Theorem~\ref{AYS-rshd}.

The next difficulty we face is that the unitary
\[
s = (s_0, s_1, \dots, s_l )
   \in \bigoplus_{k = 0}^{l} C \left( Y_k, M_{n (k)} \right)
\]
is not in the image of $A (Y)$.
(When $M$ is totally disconnected and $Y$ is both closed and open,
there is no problem: the image of $\gm_Y$ is all of
$\bigoplus_{k = 0}^{l} C \left( Y_k, M_{n (k)} \right)$.)
The cure for this problem is the following lemma, which however
requires that we look at two nested subsets $Y$ and $Z$, along with
the associated subalgebras $A (Y)$ and $A (Z)$.

\begin{lem}\label{BT-exi}
Let $M$ be a \cms\   with finite covering dimension $d$,
and let  $h \colon M \to M$ be a \mh.
Let $Y \subset M$ be closed with $\sint (Y) \neq \varnothing$.
Then every point of $\sint (Y)$ has a \nbhd\  $U \subset \sint (Y)$
such that for every closed set
$Z \subset U$ with $\sint (Z) \neq \varnothing$,
and every closed subset $S \subset \sint (Z_0)$, there is a unitary
$v \in A (Z, S)$
such that $v f = u f$ in $C^* (\Z, M, h)$ whenever
$f \in C (M)$ vanishes on $Y$.
\end{lem}

The condition on $U$ used in the proof is that there are at least
$\max \left( 1, \frac{1}{2} d \right)$ images of $\overline{U}$
under positive powers $h^r$ of $h$, with $r$ less than the smallest
first return time of $\overline{U}$ to itself, which are
contained in $\sint (Y)$.
Under this condition, the first step in the construction of $v$
is an approximate polar decomposition,
in the \rsha\   $\gm_Z (A (Z, S))$, of $u g$ for a suitable function
$g \in C (M)_{\expan_Z (S)}$ which, in particular, is required to be
equal to $1$ on $M \setminus \sint (Y)$ and to vanish on $Z$.

It isn't in general true that $\sint (Z) \neq \varnothing$ implies
$\sint (Z_0) \neq \varnothing$, although it happens that the sets we
use in the diffeomorphism case automatically have
$\sint (Z_k) \neq \varnothing$ for all $k$.

To sum up: We have what might be called the ``basic construction''
for weak approximation in $C^* (\Z, M, h)$
(not to be confused with
the basic construction of subfactor theory), namely
a triple $(Y, Z, v)$ (or a quadruple $(Y, Z, S, v)$)
consisting of closed subsets with
\[
S \subset \sint (Z_0) \subset Z \subset \sint (Y) \subset Y \subset M
\]
(or, if $S$ is not present, at least $\sint (Z) \neq \varnothing$),
and a unitary $v \in A (Z, S)$ ($A (Z)$ if $S$ is not present)
such that $v f = u f$ in $C^* (\Z, M, h)$ whenever
$f \in C (M)$ vanishes on $Y$.
We say weak approximation here because we have not approximated $u$
in norm; rather, we have a unitary $v\in A (Z, S)$ which
``acts like $u$'' (that is, like $h$) on most of the space $M$.
In particular, this construction is not the same as what we call
a ``basic approximation'' in \cite{LP2}.
The basic approximation, of which we describe an easier form in
the next section, does permit the norm approximation of $u$,
but requires two nested basic constructions and an additional unitary.

\section{An outline of the proof of local approximation}

In this section, we outline the proof of a weak form of
Theorem~\ref{main}, namely that if $h \colon M \to M$
is a minimal diffeomorphism of a
connected \csmf\  $M$ with $\dim (M) > 0$,
and if $F \subset C^* (\Z, M, h)$ is a finite set and $\ep > 0$,
then there is a \rsha\  $A \subset C^* (\Z, M, h)$
which approximately contains $F$ to within $\ep$.
This result requires most of the machinery needed for the
proof of the full direct limit decomposition result.

The crucial ingredient not yet mentioned is related to
Loring's version \cite{Lr} of Berg's technique \cite{Bg}.
This method (described in Step~7 below) requires a priori
bounds on the lengths of paths connecting certain elements in the
unitary groups of \hsa s of \rsha s.
This is an exponential length problem in the sense of \cite{Rn}.
We therefore begin by stating our exponential length result;
we require some terminology.

First, if $A$ is a unital \ca\  and $B \subset A$ is a \hsa, we define
the unitary group $U (B)$ to be
\[
U (B) = \{ u \in U (A) \colon u - 1 \in B \}.
\]
(This is the same as a common definition in terms of the unitization
$B^+$ of $B$, namely
\[
U (B) = \{ u \in U (B^+) \colon u - 1 \in B \}.
\]
Moreover, if $B$ is actually a corner,
then this group can be canonically identified
with the usual unitary group of $B$.)
Further, let
\beqr
\lefteqn{A = 
\left[ \cdots \rule{0em}{3ex} \left[ \rule{0em}{2.9ex} \left[
 {\textstyle{ C \left(X_0, M_{n (0) } \right) }}
     \oplus_{ C \left( \subrsz{ X_1^{(0)}, M_{n (1)} }
                                  \right) }
       {\textstyle{ C \left(X_1, M_{n (1)} \right) }}
                                  \right]
               \right. \right.} \\
    & &  \hspace{3em} \left. \rule{0em}{3ex} \left. \rule{0em}{2.9ex}
\oplus_{ C \left( \subrsz{ X_2^{(0)}, M_{n (2)} } \right) }
               {\textstyle{ C \left(X_2, M_{n (2)} \right) }}
       \right] \cdots \right]
     \oplus_{ C \left( \subrsz{ X_l^{(0)}, M_{n (l)} } \right) }
                {\textstyle{ C \left(X_l, M_{n (l)} \right) }}
\eeqr
be a \rsha.
If $B \subset A$ is a \hsa\  and $x \in X_k$ for some $k$, then we define
$\rank_x (B)$ to be the rank of the identity in the image of $B$ in the
finite dimensional \ca\  $M_{n (k)}$ under the map $\ev_x$ given by point
evaluation at $x \in X_k$.
If $v \in U (A)$, then we say that $\det (v) = 1$
if $\det ( \ev_x (v)) = 1$ for all $k$ and all $x \in X_k$.
(Although determinants are not well defined in \rsha s, one can show
that the condition $\det (v) = 1$ is well defined.)

\begin{thm}\label{EL-main}
Let $d, \, d' \geq 0$ be integers.
Then there is an integer $R$ such that the following holds.

Let $A$ be a \rsha\  which has a separable \rshd\  with
\tdim\  at most $d$ and \scn\  at most $d'$.
Let $B \subset A$ be a \hsa\  such that $\rank_x (B) \geq R$
for every $x$ in the total space of $A$.
Let $v \in U (B)$ satisfy $\det (v) = 1$
and be connected to $1$ by a path
$t \mapsto v_t$ in $U (B)$ such that $\det (v_t) = 1$ for all $t$.
Then there is a \ct\  path from $v$ to $1$ in $U (B)$ with length
less than $4 \pi (d' + 2)$.
\end{thm}

At this point,
we should give a brief indication of the significance of the \scn.
We explained in Section~4 of \cite{LP1} how relative versions
of the subprojection and cancellation theorems for $C (X, M_n)$
can be used to obtain analogous theorems for \rsha s.
Theorem \ref{EL-main}, however, is an exponential length theorem,
and, at a crucial step in its proof, we have only been able to prove an
approximate relative theorem for $C (X, M_n)$.
(See Theorem~6.2 of \cite{LP1}.)
Roughly speaking, errors accumulate everywhere that the \rshd\  of $A$
specifies that two algebras be glued together.
The \scn\  gives a limit on how often a \nbhd\  of a particular
point in one of the base spaces is involved in such a gluing.
It is a strengthened version of the most obvious notion
(the ``covering number''); the more obvious version proved to be
technically too weak.

The definition of the \scn\  is somewhat complicated, and is omitted;
instead, we illustrate with an example.
Let $X$ be a compact metric
space, let $E$ be a locally trivial continuous
field over $X$ with fiber $M_n$, and let $\Gm (E)$ be the
corresponding section algebra.
Then any finite cover $X_0, X_1, \dots, X_l$ of $X$
by closed subsets, such that $E |_{X_k}$ is trivial for each $k$,
induces a \rshd\  of $\Gm (E)$.
(See the proof of Proposition~1.7 of \cite{Ph6}
and Example~2.8 of \cite{LP1}.)
It can be shown that the \scn\  of this \rshd\  is the order
(as in Definition~1.6.6 of~\cite{En})
of the cover of $X$ by the sets $X_0, X_1, \dots, X_l$,
that is, the largest number $d$ such that
there are distinct $r_0, r_1, \dots, r_d$ for which
\[
\bigcap_{j = 0}^d X_{r_j} \neq \varnothing.
\]
Note the parallel with the definition of the covering dimension
(Definition 1.6.7 of \cite{En}).

At this point, we can explain how we use the condition that we have a
diffeomorphism of a \mf.
Let $Y \subset M$ satisfy $\sint (Y) \neq \varnothing$.
Our method for bounding the \scn\  requires that there be an integer
$m$ such that, for any $m + 1$ distinct integers
$r_0, r_1, \dots, r_m \in \Z$, we have
\[
\bigcap_{j = 0}^m h^{r_j} ( \partial Y ) = \varnothing.
\]
When $h$ is a minimal diffeomorphism of a compact \mf,
this is arranged as follows.
First, require that $\partial Y$ be a smooth submanifold
(of codimension $1$).
Then perturb $\partial Y$ by an arbitrarily small amount, so that
all finite sets
\[
h^{r_0} ( \partial Y ), \, h^{r_1} ( \partial Y ),
       \, \dots, h^{r_m} ( \partial Y )
\]
of distinct images of $\partial Y$ under powers of $h$
are jointly mutually transverse.
This means, first, that
$h^{r_0} ( \partial Y )$ and $h^{r_1} ( \partial Y )$
are transverse (see pages 28--30 of \cite{GP}) whenever $r_0 \neq r_1$,
so that $h^{r_0} ( \partial Y ) \cap h^{r_1} ( \partial Y )$ is a
smooth sub\mf\  (of codimension $2$;
see the theorem on page~30 of \cite{GP});
that $h^{r_2} ( \partial Y )$ and
$h^{r_0} ( \partial Y ) \cap h^{r_1} ( \partial Y )$
are transverse whenever $r_0$, $r_1$, and $r_2$ are all distinct,
so that
$h^{r_0} ( \partial Y ) \cap h^{r_1} ( \partial Y )
               \cap h^{r_2} ( \partial Y )$
is a smooth sub\mf\  (of codimension $3$); etc.
These conditions
guarantee that the intersection of any $\dim (M) + 1$ distinct
images of $\partial Y$ under powers of $h$ will be empty.
(Note, however, that the resulting upper bound on the \scn\  turns out
to be $\dim (M) ( \dim (M) + 2)$, not $\dim (M)$.
The situation is much more complicated than for
section algebras of locally trivial \ct\  fields.)
We thus have:

\begin{prp}\label{scn-t}
Let $M$ be a connected
\csmf\  with $\dim (M) = d > 0$, and let $h \colon M \to M$
be a minimal diffeomorphism.
For every $x \in M$ and open $U \subset M$ with $x \in U$,
there is a closed set $Y \subset M$ with
$x \in \sint (Y) \subset Y \subset U$
such that for every closed set $S \subset \sint (Y_0)$
(notation as in Section~2) which
is homeomorphic to a closed ball in $\R^d$,
the subalgebra $A (Y, S)$ satisfies the following properties:
\bei
\item
The \rshd\  of Theorem~\ref{AYS-rshd} has \tdim\  equal to $d$.
\item
The decomposition of Theorem~\ref{AYS-rshd}
has \scn\  at most $d (d + 2)$.
\item
In the notation of Theorem~\ref{AYS-rshd}, we have
$Y_k^{(0)} \subset \partial Y_k$ for all $k$.
\eei
\end{prp}

We hope that if $h$ is a minimal homeomorphism of a finite
dimensional compact metric space, then one might be able to
substitute a dimension theory argument for transversality in the
above.
We have not yet had time to look into this.
What to do about infinite dimensional  compact metric spaces
(such as $(S^1)^{\Z}$) is less clear.

Now we start the outline of the proof of local approximation.
We fix a connected \csmf\  $M$ with $\dim (M) > 0$ and a
minimal diffeomorphism $h \colon M \to M$.

{\bf{Step 1.}}
It suffices to prove the following:
Let
\[
f_1, f_2, \dots, f_m \in C (M) \subset C^* (\Z, M, h)
\]
be a finite collection of functions, and let $\ep > 0$.
Then there is a \rsha\  $A \subset C^* (\Z, M, h)$
which approximately contains $\{ f_1, f_2, \dots, f_m, u \}$ to within $\ep$.
(The reason is that $C (M)$ and $u$ generate $C^* (\Z, M, h)$ as a \ca.)

{\bf{Step 2.}}
Choose $\dt > 0$ so small that the functions $f_1, f_2, \dots, f_m$ are
all approximately constant to within $\frac{1}{2} \ep$ on every subset
of $M$ with diameter less than $\dt$.
Choose an integer $R$ following Theorem~\ref{EL-main} for the number
$d = \dim (M)$ and for $d' = d ( d + 2)$, and also with
$R \geq \max \left( 1, \frac{1}{2} d \right)$.
Choose an integer $N$ so large that
\[
\frac{4 \pi (d' + 2) }{N} < \ep.
\]

{\bf{Step 3.}}
Choose a quadruple $\left( Y^{(1)}, \, Z^{(1)}, \, S, \, v_1 \right)$,
as described at the end of the previous section,
consisting of closed subsets with
\[
{\textstyle{
\varnothing \neq \sint (S) \subset
 S \subset \sint \left( \rsz{Z^{(1)}_0} \right) \subset Z^{(1)}
 \subset \sint \left( Y^{(1)} \right) \subset Y^{(1)} \subset M
}}
\]
and a unitary $v_1 \in A \left( Z^{(1)}, \, S \right)$
such that $v_1 f = u f$ in $C^* (\Z, M, h)$ whenever
$f \in C (M)$ vanishes on $Y^{(1)}$.
We also require that the conclusions of Proposition~\ref{scn-t}
be satisfied.
Let $n_1 (0) < n_1 (1) < \cdots < n_1 (l_1)$
be the first return times
$n_{Y^{(1)} } (0) < n_{Y^{(1)} } (1) < \cdots
            < n_{Y^{(1)}} \left( l_{Y^{(1)}} \right)$.
We then further require that the sets involved be so small that:
\bei
\item
The sets
$Y^{(1)}, \, h^{-1} \left( Y^{(1)} \right), \, \dots,  \,
         h^{-N} \left( Y^{(1)} \right)$
are pairwise disjoint (whence $n_1 (0) > N$).
\item
The sets
$Y^{(1)}, \, h^{-1} \left( Y^{(1)} \right), \, \dots,  \,
         h^{-N} \left( Y^{(1)} \right)$
all have diameter less than $\dt$.
\item
The sets
$h (S), \, h^2 (S), \, \dots,  \, h^{n_1 (0)} (S)$
all have diameter less than $\dt$.
\item
Each of the sets
$h (S), \, h^2 (S), \, \dots,  \, h^{n_1 (0)} (S)$
is either contained in one of
$Y^{(1)}, \, h^{-1} \left( Y^{(1)} \right), \, \dots,  \,
         h^{-N} \left( Y^{(1)} \right)$
or is disjoint from all of them.
\eei
(Note that we choose $S$ after having chosen $Y^{(1)}$.)

{\bf{Step 4.}}
Choose a triple $\left( Y^{(2)}, \, Z^{(2)}, \, v_2 \right)$,
as described at the end of the previous section,
consisting of closed subsets with
\[
{\textstyle{
\varnothing \neq \sint \left( Z^{(2)} \right) \subset Z^{(2)}
 \subset \sint \left( Y^{(2)} \right) \subset Y^{(2)} \subset \sint (S)
}}
\]
and a unitary $v_2 \in A \left( Z^{(2)}, \, S \right)$
such that $v_2 f = u f$ in $C^* (\Z, M, h)$ whenever
$f \in C (M)$ vanishes on $Y^{(2)}$.
Again, we also require that the conclusions of Proposition~\ref{scn-t}
be satisfied.
Let $n_2 (0) < n_2 (1) < \cdots < n_2 (l_2)$
be the first return times
$n_{Y^{(2)}} (0) < n_{Y^{(2)}} (1) < \cdots
            < n_{Y^{(2)}} \left( l_{Y^{(2)}} \right)$.
Let $B \subset A \left( Z^{(2)} \right)$ be the \hsa\  generated by
$C_0 \left( \sint \left( Y^{(1)} \right) \right) \subset C (M)$.
We then further require that $Z^{(2)}$ be so small that
$\gm_{Z^{(2)}} (B)$, as a \hsa\  of the \rsha\  %
$\gm_{Z^{(2)}} \left( A \left( Z^{(2)} \right) \right)$,
satisfies $\rank_x \left( \gm_{Z^{(2)}} (B) \right) \geq R$
for all $x$ (in the sense discussed before Theorem~\ref{EL-main}).
(This is accomplished by requiring that
there be at least
$R$ images of $Z^{(2)}$
under positive powers $h^r$ of $h$, with $r < n_2 (0)$, which are
contained in $\sint \left( Y^{(1)} \right)$.)

{\bf{Step 5.}}
Observe that the relations
$v_j f = u f$ in $C^* (\Z, M, h)$ whenever
$f \in C (M)$ vanishes on $Y^{(j)}$ imply that
$v_1^* v_2 f = f$ whenever $f \in C (M)$ vanishes on $Y^{(1)}$.
{}From this one can deduce that $v_1^* v_2 \in U (B)$.
With the help of the condition
$\rank_x \left( \gm_{Z^{(2)}} (B) \right)
    \geq \max \left( 1, \frac{1}{2} d \right)$,
it is possible to alter the choice of $v_2$ so that, in addition to the
conditions we already have, also
$z = \gm_{Z^{(2)}} \left( v_1^* v_2 \right)
            \in U \left( \gm_{Z^{(2)}} (B) \right)$
satisfies $\det (z) = 1$
and is connected to $1$ by a path
$t \mapsto z_t$ in $U \left( \gm_{Z^{(2)}} (B) \right)$
such that $\det (z_t) = 1$ for all $t$.
(For the meaning of these conditions,
see the discussion before Theorem~\ref{EL-main}.)
Then also $v_2^* v_1 = \left( v_1^* v_2 \right)^*$
satisfies these properties.

{\bf{Step 6.}}
Apply Theorem~\ref{EL-main} to find a path
in $U (B)$ from $v_1^* v_2$ to $1$ with total length less than
$4 \pi (d' + 2)$.
Using a suitable subdivision of the domain of this path,
find unitaries
\[
v_2^* v_1 = w_0, \, w_1, \, \dots, \, w_{N - 1}, \, w_N = 1 \in U (B)
\]
such that
\[
\| w_j - w_{j - 1} \| < \frac{4 \pi (d' + 2) }{N} < \ep
\]
for $1 \leq j \leq N$.

{\bf{Step 7.}}
Define
\[
w = w_0 {\textstyle{ \left( u^{-1} w_1 u \right) }}
    {\textstyle{ \left( u^{-2} w_2 u^2 \right) }} \cdots
    {\textstyle{ \left( u^{-N} w_N u^N \right) }}.
\]
Then $w$ is a unitary in $C^* (\Z, M, h)$ with the following properties:
\bei
\item[(1)]
$w$ commutes with every $f \in C (M)$ which is constant
on each of the sets
\[
{\textstyle{
Y^{(1)}, \, h^{-1} \left( Y^{(1)} \right), \, \dots,  \,
         h^{-N} \left( Y^{(1)} \right) }}.
\]
\item[(2)]
$w$ commutes with $u v_2^*$.
\item[(3)]
$\| w v_1 w^* - v_2 \| < \ep$.
\eei
We will say something below about how these results follow.
Some of the ideas are related to calculations in Section~6 of
\cite{Pt1} and Section~2 of \cite{Pt2}.

{\bf{Step 8.}}
Set
\[
{\textstyle{
D = C^* \left( u v_2^*, \, A \left( Z^{(1)}, \, S \right) \right)
    \subset C^* (\Z, M, h)
}}
\andeqn A = w D w^*.
\]

We show that $A$ approximately contains
$f_1, f_2, \dots, f_m,$ and $u$ to within $\ep$.

Let $T_1, T_2, \dots, T_r$ consist of the sets
$Y^{(1)}, \, h^{-1} \left( Y^{(1)} \right), \, \dots,  \,
         h^{-N} \left( Y^{(1)} \right)$,
together with all of the sets
$h (S), \, h^2 (S), \, \dots,  \, h^{n_1 (0)} (S)$
which are not contained in any of
the images of $Y^{(1)}$ listed above.
By the construction in Step~(3), the sets
$T_1, T_2, \dots, T_r$ are
pairwise disjoint and have diameter less than $\dt$.
The functions $f_1, f_2, \dots, f_m$ are
all approximately constant to within $\frac{1}{2} \ep$ on every subset
of $M$ with diameter less than $\dt$ (by Step~2),
so there exist functions
$g_1, g_2, \dots, g_m \in C (M)$
which are actually constant on the sets
$T_1, T_2, \dots, T_r$
and satisfy $\| g_i - f_1 \| < \ep$ for $1 \leq i \leq m$.
These functions are then constant on all of
\[
{\textstyle{
Y^{(1)}, \, h^{-1} \left( Y^{(1)} \right), \, \dots,  \,
         h^{-N} \left( Y^{(1)} \right) }}
      \andeqn
  h (S), \, h^2 (S), \, \dots,  \, h^{n_1 (0)} (S).
\]
Now $g_i \in A \left( Z^{(1)}, \, S \right) \subset D$
and (by Step~7~(1)) $w$ commutes with $g_1, g_2, \dots, g_m$,
so $g_1, g_2, \dots, g_m \in w D w^* = A$.

We also have $w \left( u v_2^* \cdot v_1 \right) w^* \in A$.
Using the relations $w \left( u v_2^* \right) w^* = u v_2^*$
and $\left\| w v_1 w^* - v_2 \right\| < \ep$ from Step~7, we get
\[
{\textstyle{
   \left\| w \left( u v_2^* \cdot v_1 \right) w^* - u \right\|  }}
 = {\textstyle{
   \left\| w \left( u v_2^*  \right) w^* \cdot w v_1 w^*
                     - u v_2^* \cdot v_2 \right\|  }}
 < \ep.
\]
So $u$ is approximately in $A$.

{\bf{Step 9.}}
The algebra $D$, and hence $A = w D w^*$, is a \rsha\  with
\tdim\  $d$ and \scn\  at most $d' = d (d + 2)$
(that is, no more complicated than $A \left( Z^{(1)}, \, S \right)$).
This step is where $S$ is used in an essential way.

Let's assume for simplicity that $\spec \left( u v_2^* \right)$ is
the whole unit circle $S^1$.
Then it turns out that $D$ is a pullback
\[
D \cong A {\textstyle{ \left( Z^{(1)}, \, S \right) }}
 \oplus_{M_{n_1 (0)}, \ph, \ps }
         C {\textstyle{ \left( S^1, M_{n_1 (0)} \right)}}.
\]
The map
$\ps \colon
   C {\textstyle{ \left( S^1, M_{n_1 (0)} \right)}} \to M_{n_1 (0)}$
is evaluation at $1 \in S^1$.
The map
$\ph \colon
   A {\textstyle{ \left( Z^{(1)}, \, S \right) }} \to M_{n_1 (0)}$
is the evaluation on the set $S$ in the \rshd\  described in
Theorem~\ref{AYS-rshd}.
(This is really a point evaluation, because the elements of
$A {\textstyle{ \left( Z^{(1)}, \, S \right) }}$ are constant on $S$.)
The unitary $u v_2^*$ corresponds to the pair
\[
\left( 1, \, \diag \left( z, 1, \dots, 1 \right) \right)
\]
in which $z$ is the identity function $\zt \mapsto \zt$ in $C (S^1)$.

The key relation here is that $u v_2^*$ acts as $1$ off $h (S)$.
Thus, if $f \in C (M)$ vanishes on
$h (S)$, then $\left( u v_2^* \right) f = f \left( u v_2^* \right) = f$.
If in addition $f$ vanishes on $Z^{(1)}$, then
$\left( u v_2^* \right) (u f) = (u f) \left( u v_2^* \right) = u f$.
These relations imply, for example, that $u v_2^*$ commutes with
all elements of $A {\textstyle{ \left( Z^{(1)}, \, S \right) }}$.

The verification of the isomorphism with the pullback requires
lots of functional calculus.
For example, one needs to define suitable \hm s with domain
$D = C^* \left( u v_2^*, \, A \left( Z^{(1)}, \, S \right) \right)$,
or at least determine somehow all the elements of this \ca.
We omit further discussion, except to note that it is much easier
to demonstrate that there is an exact sequence
\[
0 \longrightarrow {\mathrm{Ker}} (\ph)
\longrightarrow D
\longrightarrow C {\textstyle{ \left( S^1, \, M_{n_1 (0)} \right)}}
\longrightarrow 0,
\]
as should certainly happen for a pullback with surjective maps.
This exact sequence implies (using Theorem~2.16 of~\cite{Ph6})
that $D$ is a \rsha\  with \tdim\  $d$, but doesn't
give anything about the \scn.

This finishes the outline of the proof of local approximation.

Let us now return to the explanation of Step~7.
We first explain the significance of $w$, in a greatly simplified
context---so much simplified that it does not satisfy the
hypotheses of this section.
Then we give an outline of how to prove the claimed properties
in our case.

For the simple context, let us assume that
\[
Z^{(1)} = Y^{(1)} \andeqn
M = \coprod_{j = 1}^{n} h^j {\textstyle{ \left( Z^{(1)} \right)}}.
\]
(We ignore $S$, since it is not relevant for this step.)
In this case, note that $Z^{(1)}_0 = Z^{(1)}$, that $n = n_1 (0)$,
and that $\gm_{  Z^{(1)} }$ induces an isomorphism
$A \left( Z^{(1)} \right)
           \cong M_n  \left( C \left( Z^{(1)} \right) \right)$,
under which functions constant on each of the sets
\[
{\textstyle{
Y^{(1)}, \, h^{-1} \left( Y^{(1)} \right), \, \dots,  \,
         h^{-N} \left( Y^{(1)} \right) }}
\]
are sent to the diagonal matrices in
$M_n  \left( C \left( Z^{(1)} \right) \right)$,
the last $N + 1$ diagonal entries of which are constants.
(Our simplifying assumptions imply that
$h^{-j} {\textstyle{ \left( Z^{(1)} \right)}}
          = h^{n - j} {\textstyle{ \left( Z^{(1)} \right)}}$.)

Let us further assume we have an $h$-invariant Borel probability
measure $\mu$ on $M$, and that $C^* (\Z, M, h)$ is represented
faithfully on $L^2 (M, \mu)$ with $C (M)$ acting as multiplication
operators and $u$ acting as $u \xi = \xi \circ h^{-1}$.
There is a direct sum decomposition
\[
L^2 (M, \mu) = \bigoplus_{j = 1}^{n}
    L^2 {\textstyle{ \left( h^j \left( Z^{(1)} \right) \right)}},
\]
which determines an identification of $L ( L^2 (M, \mu) )$ with
$M_n {\textstyle{ \left( L^2 \left(  Z^{(1)} \right) \right)}}$
which is compatible in a suitable sense with the isomorphism
$\gm_{  Z^{(1)} }$.
Further let $e_j$ be the projection onto
$L^2 {\textstyle{ \left( h^j \left( Z^{(1)} \right) \right)}}$.
With respect to this identification, we can write
\[
u = \left( \begin{array}{ccccccc}
         0 &      0 & \cdots & \cdots &      0 &      0 &  u^{(0)} \\
         1 &      0 & \cdots & \cdots &      0 &      0 &      0 \\
         0 &      1 & \cdots & \cdots &      0 &      0 &      0 \\
    \vdots & \vdots & \ddots &        & \vdots & \vdots & \vdots \\
    \vdots & \vdots &        & \ddots & \vdots & \vdots & \vdots \\
         0 &      0 & \cdots & \cdots &      1 &      0 &      0 \\
         0 &      0 & \cdots & \cdots &      0 &      1 &      0
 \end{array} \right),
\]
with $u^{(0)} \in e_1 L ( L^2 (M, \mu) ) e_n$.
(Note that it is equal to the shift matrix $s_0$ considered in
Section~2, except for the upper right corner.)
Similarly, we can write
\[
v_j = \left( \begin{array}{ccccccc}
         0 &      0 & \cdots & \cdots &      0 &      0 &  v_j^{(0)} \\
         1 &      0 & \cdots & \cdots &      0 &      0 &      0 \\
         0 &      1 & \cdots & \cdots &      0 &      0 &      0 \\
    \vdots & \vdots & \ddots &        & \vdots & \vdots & \vdots \\
    \vdots & \vdots &        & \ddots & \vdots & \vdots & \vdots \\
         0 &      0 & \cdots & \cdots &      1 &      0 &      0 \\
         0 &      0 & \cdots & \cdots &      0 &      1 &      0
 \end{array} \right).
\]
Again, the difference is in the the upper right corner,
but note that $v_1$ and $v_2$ are now in $A \left( Z^{(2)} \right)$.

In this situation, we let $w'_i = e_0 w_i e_0$, and identify $w$ as
\[
w = \diag \left( 1, \, 1, \, \dots, \, 1, \, w'_N, \, w'_{N - 1}, \,
    \dots, \, w'_1, \, w'_0 \right).
\]
(We have used the fact that $n \geq N + 1$.)
Now Condition~(1) of Step~7
follows from the fact that $w$ is block diagonal
and that functions in $C (M)$ constant on each of the sets
\[
{\textstyle{
Y^{(1)}, \, h^{-1} \left( Y^{(1)} \right), \, \dots, \,
         h^{-N} \left( Y^{(1)} \right) }}
\]
are diagonal matrices,
the last $N + 1$ diagonal entries of which are constants.
For Condition~(2) of Step~7, we calculate:
\[
u v_2^*
   = \diag \left( u^{(0)}
                {\textstyle{ \left( \rsz{ v_2^{(0)} } \right)^*}}, \,
 1, \, 1, \, \dots, \, 1 \right).
\]
This element clearly commutes with $w$.
(The worst case is $n = N + 1$; then, recall that $w_N = 1$.)
For Condition~(3) of Step~7, we estimate instead
$\left\| w - v_2 w v_1^* \right\|$.
(This is easily seen to be equivalent.)
A computation shows that
\begin{align*}
v_2 w v_1^*
  & = \diag \left( v_2^{(0)}
            w'_0 {\textstyle{ \left( \rsz{ v_1^{(0)} } \right)^*}},
    \, 1, \, \dots, \, 1, \, 1, \, w'_N, \,
    \dots, \, w'_2, \, w'_1 \right)  \\
  & = \diag \left( \rsz{ 1, \, 1, \, \dots, \, 1, \, 1, \, w'_N, \,
    \dots, \, w'_2, \, w'_1 } \right).
\end{align*}
(The entries of $w$ have all been moved one space down the
diagonal.
In addition, the new first entry has been modified.
Since $w_0 = v_2^* v_1$, we have
$w'_0 = {\textstyle{\left( \rsz{ v_2^{(0)} } \right)^* v_1^{(0)} }}$.)
Therefore, using $w_N = 1$, we get
\[
{\textstyle{ \left\| \rsz{w - v_2 w v_1^*} \right\| }}
  = \max_{1 \leq j \leq N}
         {\textstyle{ \left\| w_j -  w_{j - 1} \right\| }}
 < \ep.
\]

In the actual situation, we work inside $C^* (\Z, M, h)$.
Let $B \subset C^* (\Z, M, h)$ be the \hsa\  of Step~4.
For the matrix decomposition, we substitute the fact that the
\hsa s
\[
B, \, u^{-1} B u, \, u^{-2} B u^2, \, \dots, \, u^{-N} B u^N
\]
are orthogonal in $C^* (\Z, M, h)$.
This follows from the fact that the sets
\[
{\textstyle{
Y^{(1)}, \, h^{-1} \left( Y^{(1)} \right), \, \dots,  \,
         h^{-N} \left( Y^{(1)} \right) }}
\]
are pairwise disjoint.
As a consequence, the factors
\[
w_0, \, u^{-1} w_1 u, \, u^{-2} w_2 u^2, \, \dots, \,
    u^{-N} w_N u^N
\]
of $w$, which are in the unitary groups of these \hsa s,
all commute with each other, and also with any function
$f \in C (M)$ which is constant on each of the sets
\[
{\textstyle{
Y^{(1)}, \, h^{-1} \left( Y^{(1)} \right), \, \dots,  \,
         h^{-N} \left( Y^{(1)} \right) }}.
\]
When proving that $w$ commutes with $u v_2^*$,
it helps to show first that
\[
u^{-j} w_j u^j = v_2^{-j} w_j v_2^j
\]
for $0 \leq j \leq N$.
In fact, this is true if $w_j$ is replaced by any $b \in C^* (\Z, M, h)$
which differs by a scalar from an element of $B$.
For the verification of the norm estimate in Condition~(3) of Step~7,
one needs in addition the following fact, which is the analog of the
estimate on the difference of diagonal matrices above:
if $C_0, C_1, \dots, C_N$ are orthogonal \hsa s in a \ca\  $A$,
and if $y_j, \, z_j \in U (C_j)$ for $0 \leq j \leq N$,
then
\[
\| y_0 y_1 \cdots y_N - z_0 z_1 \cdots z_N \|
   =  \max_{0 \leq j \leq N}
         {\textstyle{ \left\| y_j -  z_j \right\| }}.
\]

\section{Direct limit decomposition}

We give here a very brief approximate outline of the modifications
necessary to achieve the direct limit decomposition of Theorem~\ref{main},
as opposed to merely local approximation.
The previous section describes the construction of a
(simple version of) a single ``basic approximation'',
and the problem is to arrange successively better ones so as to obtain
an increasing sequence of subalgebras of $C^* (\Z, M, h)$.
As will be clear, putting everything together requires complicated
notation, and there are interactions between the modifications
described below which we do not have room to discuss here.

First, the unitary corresponding to $w$ in each new
basic approximation must commute with all elements of
the subalgebra $A \left( Z^{(2)} \right)$ from the previous one.
This requires two changes.
The old subalgebra $A \left( Z^{(2)} \right)$ must be replaced
by $A \left( Z^{(2)}, \, T \right)$ for some suitable $T$,
and the new set $Y^{(1)}$ must be contained in $T$.
Also, the sequence
\[
v_2^* v_1 = w_0, \, w_1, \, \dots, \, w_{N - 1}, \, w_N = 1 \in U (B)
\]
used to construct the new $w$ must now consist of constant subsequences,
the lengths of which are certain return times associated with the old
$Z^{(2)}$.

Second, having constructed one approximating subalgebra, say $A_0$,
the next one, say $A_1$, will be slightly ``twisted'' with respect
to $A_0$, even with the adjustment above.
To straighten this out, it is necessary to modify $A_0$ by
replacing $v_2$ in the construction by a nearby unitary.
Then, after constructing $A_2$, one must further modify the
unitaries $v_2$ associated with both $A_1$ and $A_0$, etc.
Enough control must be maintained that the sequences of modifications
converge to unitaries not too far from the original choices.

Third, even apart from the ``twisting'' referred to in the
previous paragraph, the use of the subsets $S$ leads to problems
with the expected inclusion relations between subalgebras.
Suppose, for example,
we have closed subsets $Y$ and $Z$, satisfying the conclusions
of Proposition~\ref{scn-t}, with associated first return times
\[
n_Y (0) < n_Y (1) < \cdots < n_Y (l_Y) \andeqn
n_Z (0) < n_Z (1) < \cdots < n_Z (l_Z),
\]
and with corresponding subsets
\[
Y_0, \, Y_1, \, \dots, \, Y_{l_Y} \subset Y \andeqn
Z_0, \, Z_1, \, \dots, \, Z_{l_Z} \subset Z.
\]
Suppose that
\[
\varnothing \neq S \subset \sint (Z_0) \subset Z \subset \sint (Y_0)
\]
(in particular, $Z \subset Y$), and that $n_Z (0) > n_Y (0)$
(this is the relevant situation, because arbitrarily good
approximations require arbitrarily large values of the
smallest first return time).
We have $A (Y) \subset A (Z)$, because every function in $C (M)$
which vanishes on $Y$ also vanishes on $Z$.
However, it is not true that $A (Y, S) \subset A (Z, S)$.
In fact, $C (M) \cap A (Z, S)$ consists of those functions in $C (M)$
that are constant on the sets
\[
h (S), \, h^2 (S), \, \dots, \, h^{n_Z (0)} (S),
\]
$C (M) \cap A (Y, S)$ consists of those functions in $C (M)$
that are constant on the sets
\[
h (S), \, h^2 (S), \, \dots, \, h^{n_Y (0)} (S),
\]
and $n_Y (0) < n_Z (0)$, so
$C (M) \cap A (Y, S) \subsetneqq C (M) \cap A (Z, S)$.

To fix this problem, it is necessary to replace the single set
$S$ in the construction of $A (Y, S)$ by a whole family of subsets.
One must require that whenever $h^j (S) \subset Y$, with
$0 < j < n_Z (0)$, then there is $k$ with $h^j (S) \subset \sint (Y_k)$.
Then one uses the collection of all such $h^j (S)$, rather than just
$S$ itself, with the obvious modification to account for the fact
that they are no longer all subsets of $\sint (Y_0)$.
The resulting subalgebra is a proper subalgebra of $A (Y, S)$.

In the inductive construction of an increasing sequence of
approximating subalgebras of $C^* (\Z, M, h)$, this works out as
follows.
First, one constructs an approximating algebra $A_0^{(0)}$.
Then one constructs an approximating algebra $A_1^{(1)}$,
incorporating the first two modifications
discussed above, and using a sufficiently small set $S$.
Next, one replaces $A_0^{(0)}$ by a smaller algebra $A_0^{(1)}$,
using the approach
outlined in the previous paragraph on the algebra
$A \left( Z^{(1)}, \, S \right)$ appearing in the definition
of $A_0^{(0)}$,
but with the set $S$ from the construction of $A_1^{(1)}$.
That done, one constructs $A_2^{(2)}$.
Then it is necessary to go back and replace both
$A_1^{(1)}$ and $A_0^{(1)}$ (in that order) by smaller subalgebras
$A_1^{(2)}$ and $A_0^{(2)}$,
in a similar way.
This procedure continues for all $n$.

There are two problems.
First, $\bigcap_{k = n}^{\infty} A_n^{(k)}$ must still be large enough
to approximate not too badly the finite set that the first algebra
$A_n^{(n)}$ was constructed to approximate.
Second, $\bigcap_{k = n}^{\infty} A_n^{(k)}$ must still be a
\rsha\  with \tdim\  at most $d$ and \scn\  at most $d (d + 2)$.
Since subalgebras of \rsha s need not even be \rsha s
(see Example~3.6 of \cite{Ph6}), this requires work.
The construction of the subalgebra $A (Y, S)$ can be viewed
as identifying the subset $S$ of $Y$ to a point.
By the time the inductive process of the previous paragraph is
complete, one must identify infinitely many subsets of $Y$ to
(distinct) points,
in such a way that the resulting space is not only Hausdorff
(there is trouble even here) but in fact has dimension no
greater than $\dim (Y)$.
The details are quite messy.

\end{document}